\def\Date{01.05.2011}

\magnification    = \magstep1
\input amstex
\input epsf
\input texdraw
\documentstyle{amsppt}
\vsize=20.8cm \hsize=14.2cm \hoffset=-0.4cm \voffset=-0.5cm
\emergencystretch = 10 pt
\hfuzz        =  8 pt
\TagsOnRight
\NoBlackBoxes

\def\pt(#1,#2){\move(#1 #2) \fcir f:0 r:\ptrad}
\def\ptsize(#1){\def\ptrad{#1}} \ptsize(0.25)

\def\dl(#1,#2,#3,#4){\move(#1 #2) \rlvec(#3 #4)} 
\def\vku(#1,#2,#3,#4,#5){
  \move(#1 #2) \rlvec(#3   0) \rlvec(#4  #4) \rlvec( #4 -#4) \rlvec( #5 0)}
\def\vkd(#1,#2,#3,#4,#5){
  \move(#1 #2) \rlvec(#3   0) \rlvec(#4 -#4) \rlvec( #4  #4) \rlvec( #5 0)}
\def\dmd(#1,#2,#3){
  \move(#1 #2) \rlvec(#3  #3) \rlvec(#3 -#3) \rlvec(-#3 -#3) \rlvec(-#3 #3)}
\def\rct(#1,#2,#3,#4){
  \move(#1 #2) \rlvec(#3   0) \rlvec( 0 -#4) \rlvec(-#3   0) \rlvec(  0 #4)}
\def\tru(#1,#2,#3,#4,#5){
  \move(#1 #2) \rlvec(#3  #4) \rlvec( 0 -#5) \rlvec(-#3 -#4) \rlvec(  0 #5)}
\def\trd(#1,#2,#3,#4,#5){
  \move(#1 #2) \rlvec(#3 -#4) \rlvec( 0 -#5) \rlvec(-#3  #4) \rlvec(  0 #5)}


\everytexdraw{\drawdim mm \arrowheadtype t:V}

\def\a{\alpha}          \def\b{\beta}        \def\({\left(}
\def\g{\gamma}          \def\G{\Gamma}       \def\){\right)}
\def\d{\delta}                 \def\[{\left[}
\def\la{\lambda}             \def\]{\right]}
               
\def\p{\partial}        \def\Om{\Omega}

 \def\Sum{\Sigma}
\def\maxo{\vee}         \def\mino{\wedge}

\def\F{{\Cal F}}                 
\def\W{{\Cal W}}                 
\def\Re{{\Bbb R}}                
\def\Rn{{\Re^n}}                 

\def\card{\operatorname{card}}   
\def\diam{\operatorname{diam}}   
\def\len{\operatorname{len}}     

\def\mcig/{mCigar}               
\def\wsli/{wslice}               
\def\wslid/{wslice domain}       
\def\Wsli/{wslice${}^+$}         
\def\dwa{\text{WS}_\a}           

\def\cl[#1]{{\overline{#1}}}     
\def\ls#1,#2;{[#1\to #2]}        

\def\approxx#1{\;\raise.4ex\hbox{$#1$}\mskip-14mu\lower.7ex%
   \hbox{$\sim$}\;}              

\def\lc{{ \approxx < }}          

\def\QED{\ifmmode\eqno\qed\else\hfill\ \qed\fi}

\def\open#1;{{\par\noindent\bf #1.\, }}

\comment
+--------------  Customizable horizontal line  --------------+
| "\hr f,e,w,b": f=skip fwd, e=elevate, w=width, b=skip back |
+------------------------------------------------------------+
\endcomment
\def\hr#1,#2,#3,#4:%
{\mskip #1 mu \overline{\raise #2 ex \hbox to #3 pt{}} \mskip -#4 mu}

\comment
+---------  New "average"-symbol (line-through-integral)  ---------+
|                       Options in form:                           |
|       Display                                 Text               |
|       Script                                  Scriptscript       |
+------------------------------------------------------------------+
\endcomment
\def\av_#1{\mathchoice
{\hr 8.0,0.2,5.0,16.3: \int_{#1}}       {\hr 6.0,0.2,3.0,11.6: \int_{#1}}
{\hr 5.0,0.2,2.5,11.3: \int_{#1}}       {\hr 5.0,0.1,2.5,12.5: \int_{#1}}
}

\def\Item#1;{\item"(#1)"} 

\def\ds{\displaystyle}

\document
\topmatter
\title
Distinguishing properties of weak slice conditions II
\endtitle
\author  Stephen M. Buckley, Andr\'e Diatta, and Alexander Stanoyevitch
\endauthor
\address
Department of Mathematics and Statistics, National University of Ireland
Maynooth, Maynooth, Co. Kildare, Ireland.          \endaddress%
\email sbuckley\@maths.may.ie    \endemail
\address
Department of Mathematical Sciences, Mathematics and Oceanography Building,
Peach Street, Liverpool. L69 7ZL, United Kingdom   \endaddress%
\email adiatta\@liv.ac.uk        \endemail
\address
Department of Mathematics, California State University-Dominguez Hills, 1000
E. Victoria Street, Carson, CA  90747, USA
\endaddress%
\email astanoyevitch\@csudh.edu   \endemail%
\abstract{ The slice condition and the more general weak slice conditions
are geometric conditions on Euclidean space domains which have evolved over
the last several years as a tool in various areas of analysis. This paper
examines some of their finer distinctive properties. }
\endabstract
\thanks The first and second authors were partially supported by Enterprise
Ireland. \hfill\break\indent \Date             \endthanks
\endtopmatter

\parskip      =  4 pt plus 2 pt minus 2 pt
\parindent    = 15 pt
\baselineskip = 22 pt plus 0.5 pt minus 1.0 pt  
\baselineskip = 12 pt plus 1.0 pt minus 0.1 pt  

\def\intrsec{0} \def\termsec{1} \def\wslisec{2}
\def\ZvsPsec{3} \def\AvsBsec{4}

\head\intrsec. Introduction \endhead

The {\it slice condition} is a metric-geometric condition for domains in
Euclidean spaces $\Rn$.  It is a very weak condition which, in particular,
is satisfied by every simply connected planar domain, and was introduced by
the first author and Koskela \cite{BK2} to obtain a set of geometric
classifications of domains in Euclidean spaces which support any of the
Sobolev imbeddings, $p\ge n$. In later research, variations of the slice
condition, including the weaker conditions known as {\it weak slice
conditions} were used to refine these results and also to investigate
questions in other areas of analysis; see \cite{BO}, \cite{BS1},
\cite{BS2}, \cite{B1}, \cite{B2}, \cite{BB}. In particular, it is shown in
\cite{BB} that in many metric measure spaces, including Euclidean space,
one version of the slice condition is equivalent to Gromov hyperbolicity.
This version implies all other slice-type conditions in the literature, so
we may think of all slice-type conditions as weak versions of Gromov
hyperbolicity.

With this range of applications, it should be useful to have a solid
understanding of (weak) slice conditions, and in particular whether and how
they differ from one another. Many properties and examples of these
conditions were obtained in \cite{BS1} and \cite{BS2} but some fundamental
questions remained, including a few that were listed in Section~6 of
\cite{BS2} as open problems. A couple of these questions were answered in
\cite{BS3}. In this paper, we construct examples to answer two of the
remaining open problems in \cite{BS2}.

After some basics in Section~\termsec, we define and briefly discuss the
weak slice conditions in Section~\wslisec. Our first example is given in
Section~\ZvsPsec: it shows that there are $0$-\wsli/ domains (i.e. weak
slice domains with a certain parameter $\a$ equal to zero), which are not
slice domains, resolving Open Problem~C in \cite{BS2}. In Section~\AvsBsec,
we show that for any pair of distinct numbers $\a,\b\in[0,1)$, there is a
domain which is an $\a$-\wsli/ domain but not a $\b$-\wsli/ domain, thereby
resolving Open Problem~B in \cite{BS2}. When $\a\ge\b$, this is not hard to
deduce from the results in \cite{BS2}, but it is somewhat surprising that
the same is true when $\a<\b$. In fact we prove the following result.

{\bf Theorem A.} {\it For each $0<\a_0<1$, there are bounded Euclidean
domains $\Om_1$ and $\Om_2$ such that $\Om_i$ is an $\a$-\wsli/ domain,
$0\le\a<1$, if and only if $\a\le\a_0$ (if $i=1$) or $\a\ge\a_0$ (if $i=2$).
}

\head\termsec. Notation and Terminology \endhead

Throughout this paper we will consistently employ the following notation.
Note that certain parameters are optional in the sense that they are
omitted from the notation when understood or when the exact choice is
unimportant.

$(\Om,d)$ is a rectifiably connected incomplete metric space possibly
subject to additional restrictions (it is often just a domain in Euclidean
space), $\cl[\Om]$ is its metric completion (viewed as a superset of
$\Om$), and $\p\Om:=\cl[\Om]\setminus\Om$. For points $x,y\in\Om$, a set
$E\subset\Om$, positive numbers $r,s$, we let:

$r\maxo s$ and $r\mino s$ denote the maximum and minimum, respectively, of
$r$ and $s$;

$\lceil r \rceil$ and $\lfloor r \rfloor$ denote the smallest integer $m\ge
r$, and the largest integer $m\le r$, respectively;

$\len(E)\equiv \len_d(E)$ denotes the Hausdorff $1$-dimensional measure of
$E$ with respect to the metric $d$ (so if $E$ is an arc, $\len_d(E)$ is
just its $d$-arclength);

$\diam(E)\equiv\diam_d(E)$ denotes the $d$-diameter of $E$;

$\d(x)\equiv\d_\Om(x)$ denotes the distance from $x$ to $\p\Om$,

$B(x,r)\equiv B_{d,\Om}(x,r):=\{y\in\Om \,:\, d(x,y)<r\}$,

$B_x:=B_d(x,\d_\Om(x))$, and

$\G_\Om(x,y)$ denotes the class of all rectifiable paths $\la:[0,t]\to\Om$
for which $\la(0)=x$ and $\la(t)=y$. We do not distinguish notationally
between paths and their images. Whenever $E$ is an (open or closed) ball,
$tE$ denotes its concentric dilate by a factor $t>0$.

For $\a\in[0,1]$ we will also make extensive use of {\it
subhyperbolic lengths} and the corresponding {\it metrics}. Denoting
arclength measure by $ds$, we define these quantities by
$$
\align%
\len_\a(\g)\equiv\len_{\a,\Om}(\g) &:= \int_\g \d_{\Om}^{\a-1}(z)\,ds(z),
\qquad\text{whenever $\g$ is a rectifiable path in $\Om$} \\
d_{\a,{\Om}}(x,y)  &:=\inf_{\g\in\G_{\Om}(x,y)}\,\len_{\a,\Om} (\g),
\endalign
$$
We note that if $\Om$ is a domain in Euclidean space, or in an imbedded
$k$-manifold in $\Rn$, then $\len_{0,\Om}$ and $d_{0,\Om}$ are the
well-known {\it quasihyperbolic length} and {\it quasihyperbolic distance},
and $d_{1,\Om}$ is the {\it inner metric} with respect to $\Om$. For
brevity, we shall denote the inner metric on $\Om$ as $d_{\Om}$ and the
corresponding inner diameter of a subset $E$ of $\Om$ as $\diam_\Om(E)$ in
such cases. We shall also write $k(x,y)\equiv k_\Om(x,y)$ in place of
$d_{0,\Om}(x,y)$.

Let us call $\g\in\G_\Om(x,y)$ {\it $(\a;C_1,C_2)$-efficient}, or simply
{\it $\a$-efficient}, if
$$ \len_{\a,\Om}(\g)\le (1+C_1)d_{\a,\Om}(x,y)+C_2\,$$
We say that $\g\in\G_\Om(x,y)$ is an {\it $(\a,C_1,C_2)$-quasigeodesic} for
$x,y$ if $\g$ and all its subpaths are $(\a;C_1,C_2)$-efficient, while we
say that $\g$ is an {\it $\a$-geodesic} if it is $(\a;0,0)$-efficient (or
equivalently an $(\a;0,0)$-quasigeodesic). Obviously, efficient paths always
exist, with $(C_1,C_2)$ as close to $(0,0)$ as we wish, but $\a$-geodesics
might not exist. For instance in the Euclidean case, $\a$-geodesics exist if
$\a=0$, but might not if $\a>0$; see \cite{GO} and \cite{BS1, Example~1.1}.

Let $C\ge 1$, $x,y\in\Om$, and let $\g\in\G_\Om(x,y)$ be a path of length
$l$ which is parametrized by arclength. We say that $\g$ is a {\it
$C$-uniform path} for $x,y\in\Om$ if $\,l\le Cd(x,y)$ ({\it bounded turning
condition}) and $\,t\mino(l-t)\le C\d_\Om(\g(t))$ ({\it cigar condition}).
In this case, we get the following estimates
\def\unifest{\termsec.1}
$$
d_{\a,\Om}(x,y) \le \cases
  4C^2\log\(1+\ds{d(x,y)\over \d_\Om(x)\mino\d_\Om(y)}\), &\a=0, \\
  C'[\d_\Om(x)\maxo\d_\Om(y)\maxo d(x,y)]^\a, &0<\a\le 1.
\endcases
\tag\unifest
$$
where $C'=C'(C,\a)$. The $\a>0$ case follows by an easy integration,
estimating distance to the boundary by the triangle inequality for the
initial and final parts of the path that are close to $x$ and $y$,
respectively, and by uniformity for the rest of the path. The case $\a=0$ is
Lemma~2.14 of \cite{BHK}.



\head\wslisec. Weak Slice and Slice Conditions \endhead

In this section we define, and briefly discuss, weak slice conditions;
throughout we assume that $0\le\a<1$. For more details, we refer the reader
to \cite{BS1}, \cite{BS2}, and \cite{BS3}. We also define the slice
condition.

Suppose $C\ge 1$. A finite collection $\F$ of pairwise disjoint open subsets
of $\Om$ is a {\it set of $C$-wslices} for $x,y\in\Om$ if
$$
\align
&\forall\;S\in\F\;\;\forall\;\la\in\G_\Om(x,y):\quad
  \len(\la\cap S)\ge d_S/C \tag WS-1 \\
&\forall\;S\in\F:\;
  S\cap B(x,\d(x)/C) = S\cap B(y,\d(y)/C) = \emptyset,\tag WS-2
\endalign
$$
where $d_S\ge\diam(S)$ is some finite number associated with each wslice
$S$. We refer to such a set of data $\{(S,d_S)\mid S\in\F\}$ as being {\it
$C$-admissible} for the pair $x,y\in\Om$. Next, we define $\dwa(x,y;\Om;C)$
by
$$
\align
\dwa(x,y;\Om;C):=\sup\{\,&\d_\Om^\a(x)+\d_\Om^\a(y)+\sum_{S\in\F} d_S^\a: \\
&\quad\{(S,d_S)\mid S\in\F\}\text{ is $C$-admissible for } x,y\in\Om\,\}
\endalign
$$
%
%
Note that $\dwa(x,y;\Om;C)\ge\d_\Om^\a(x)+\d_\Om^\a(y)$, since the empty
set is trivially $C$-admissible.  {\it A priori}, $\dwa(x,y;\Om;C)$ could
possibly be infinite, but, at least in the Euclidean context, it is
bounded. In fact, Lemma~2.3 of \cite{BS1} implies that there exists a
constant $C'=C'(C,\a)$ such that
$$ \dwa(x,y;\Om;C)\le C'[\d_\Om^\a(x)+\d_\Om^\a(y)+d_{\a,\Om}(x,y)]. $$
We use subscript notation such as $\F:=\{S_i\}_{i=1}^m$ and $d_i:=d_{S_i}$ in
cases where we know that $\F$ is nonempty.

We define an $\a$-\wsli/ space essentially by reversing this last inequality
for large subhyperbolic distance. More precisely, we say that the pair $x,y$
satisfy an {\it $(\a,C)$-\wsli/ condition}, $C\ge 1$, if
$$ d_{\a,\Om}(x,y)\;\le C\,\dwa(x,y;\Om;C), \tag WS-3 $$
and we say that $\Om$ is a (two-sided) {\it $(\a,C)$-\wsli/ space} if all
pairs of points in $\Om$ satisfy an $(\a,C)$-\wsli/ condition\footnote{ In
\cite{BS1} and \cite{BS2}, the labels (WS-2) and (WS-3) were reversed, but
that does not suit our more general discussion here.}. When $\a=0$, (WS-3)
simply says that $k(x,y)\le C(2+\card(\F))$, where $\F$ is a $C$-wslice
collection of maximal cardinality.  Note that in light of (WS-1),
each of the slices $S$ must separate $x$ from $y$ in $\Om$. It is also
convenient to say that a $C$-admissible set $\{(S,d_S)\mid S\in\F\}$ for
$x,y\in\Om$ is an {\it $(\a,C)$-\wsli/ dataset} for $x,y$ if we
additionally have the following condition:
$$
d_{\a,\Om}(x,y)\le C\,\(\d_\Om^\a(x)+\d_\Om^\a(y)+\sum_{S\in\F} d_S^\a \)
$$
If the numbers $d_S$ are not specified, it is assumed that
$d_S:=\diam_d(S)$.

Oftentimes the value of the constant $C$ is unimportant and so we will on
such occasions refer simply to ``$\a$-\wsli/ conditions and/or
domains". Modulo a possible augmentation of $C$, condition (WS-2) can
actually be dropped in case $\a>0$, but it is essential in case $\a=0$, lest
every domain be a $(0,C)$-\wsli/ domain; see \cite{BS2, Theorem~5.1}.

In working with the weak slice conditions, the following additional
hypotheses have often turned out to be useful:
$$
\alignat 3
&\forall\;S\in\F\;\;\exists\;(\a;C_1,0)\text{-efficient }\g\in\G_\Om(x,y):
\quad&&\len_{\a,\Om}(\g\cap S)\le Cd_S^\a\tag WS-4 \\
&\forall\;S\in\F\;\;\exists\;z_S\in S:\quad &&B_d(z_S,d_S/C)\subset S
\tag WS-5 \\
&\forall\;S\in\F\;\;\forall\;\la\in\G_\Om(x,y):\quad&&\diam_d(\la_S)\ge
d_S/C, \tag WS-$1^+$
\endalignat
$$
where $\la_S$ denotes a component of $\la\cap S$ of maximal diameter. We
refer to $(\a,C)$-\wsli/ domains which satisfy (WS-4), (WS-5), and
(WS-$1^+$) as $(\a,C)$-\Wsli/ domains. Of these extra conditions, only
(WS-$1^+$) is significant if we do not care about the exact value of $C$,
since, modulo a possible quantitative change in the value of $C$, (WS-4) and
(WS-5) can be assumed without loss of generality; see \cite{BS3, Section
2}. The choice of $C_1>0$ and $\g$ in (WS-4) is unimportant; we can even
take $\g$ to be an $\a$-geodesic (and so $C_1=0$) if one exists. We suspect
(at least in a Euclidean or inner Euclidean context, and modulo a controlled
increase in the value of $C$ and a change in the \wsli/ dataset) that
(WS-$1^+$) also follows from the $(\a,C)$-\wsli/ condition, but we cannot
prove this.

If $\Om\subsetneq\Rn$ is a domain, we call $\Om$ an $(\a,C)$-\wsli/, or
inner $(\a,C)$-\wsli/, domain if it is an $(\a,C)$-\wsli/ space with respect
to the Euclidean or inner Euclidean metric, respectively. Notice that the
difference between Euclidean and inner Euclidean $\a$-\wsli/ domains is
rather minor since distance to the boundary, the associated subhyperbolic
metrics and the Hausdorff $1$-dimensional measure are unchanged, and so
there is no difference in any of (WS-1) through (WS-5). The only change is
in the requisite lower bound in the size of the $d_S$ (from $\diam(S)$ to
$\diam_\Om(S)$). Nevertheless, Example~3.1 of \cite{BS3} shows that there
are \wsli/ domains that are not inner \wsli/ domains.


%

We say that the pair $x,y\in\Om$ satisfy the {\it $C$-slice condition},
$C\ge 1$, if there exists $\F:=\{(S_i,d_i)\}_{i=1}^m$, with
$d_i\equiv\diam_d(S_i)$, and an $(0;C-1,0)$-efficient path
$\g\in\G_\Om(x,y)$ such that:%
\roster%
\Item a; $\F$ is an $(\a,C)$-\wsli/ dataset for $x,y$;%
\Item b; (WS-4) and (WS-5) hold for each $1\le i\le m$, $\a=0$;%
\Item c; $\forall\;1\le i\le m,\; z\in\g\cap S_i:\quad 1/C\le
  \d_\Om(z)/d_i\le C$;\smallskip%
\Item d; $\g\subset B_{k_\Om}(x,C)\cup B_{k_\Om}(y,C)\cup
  \(\bigcup_{i=1}^m\cl[S_i]\)$.
\endroster%
Slice spaces and domains are then defined in the same manner as their weak
slice equivalents.

This definition of a slice condition is different from the original (inner)
Euclidean definition in \cite{BK2}, but is equivalent to it in the Euclidean
and inner Euclidean settings (modulo a quantitatively controlled change in
$C$). For the interested reader, we note that the original definition
implies (a) by \cite{BS1, Lemma~2.4}, while (b)--(d) are easy to deduce from
the original definition. In the original definition, the path $\g$ is not
assumed to be $0$-efficient, but this follows from the previously mentioned
estimate $\dwa(x,y;\Om;C)\lc \d_\Om^\a(x)+\d_\Om^\a(y)+d_{\a,\Om}(x,y)+1$.
Proving that the original definition follows from the new one is routine.
We point out that we still do not know whether (WS-$1^+$) holds for slice spaces
(see Open Problem A in Section 6 of \cite{BS2}).

In the (inner) Euclidean setting, we point out that the (important) upper
bound of (c) is redundant. Indeed, Lemma 2.2 of \cite{BS1} tells us that if
$\{S_i,d_i\}_{i=1}^m$ is a $(0,C)$-\wsli/ dataset for points $x,y$ in a
Euclidean domain $\Om$, and $d$ is either the Euclidean or inner Euclidean
metric, then $\d_\Om(w) < C\diam_d(S_i)$ for all $w\in S_i$, $1\le i\le m$.

The point of our new definition is that it emphasizes the distinction
between slice and $0$-\wsli/ conditions. Since (WS-4) and (WS-5) follow
quantitatively from any $\a$-\wsli/ condition, it seems that the crucial
distinction is the existence of a path $\g$ which is covered by the closure
of the slices and quasihyperbolic balls around $x,y$. Intuitively, this
means that we are able to ``slice up nicely all of the region between $x$
and $y$'', whereas in a $0$-\wsli/ condition, we merely assume that we can
``slice up nicely a reasonably large part of the region between $x$ and
$y$''.

\head\ZvsPsec. $0$-\wsli/ but not slice \endhead

\def\figA{Figure~\ZvsPsec.1}
\def\exaA{Example~\ZvsPsec.2}

Here we give an example of a $0$-\wsli/ domain that is not a slice domain,
thereby resolving Open Problem~C in \cite{BS2}. Simpler examples with
related properties can be found elsewhere. Specifically, Proposition 4.5 of
\cite{BS1} allows one to construct examples of $\a$-\wsli/ domains, $\a>0$,
that are not slice domains; in fact, they are not even $0$-\wsli/ domains.
A one-sided $0$-\wsli/ domain (meaning that (WS-3) is assumed for arbitrary
$x$ and a fixed $y$) that is not a one-sided slice domain is given in
\cite{B1, Example~4.9}. However, examples similar to these cannot lead to a
(two-sided) $0$-\wsli/ domain that is not a slice domain.


\bigskip\medskip
\centertexdraw{
  \setunitscale 0.68

  \linewd 0.15 \lpatt(0.5 1.2)
  \rct( 36,66,3,12)  \ifill f:0.7 
  \tru( 60,76,3,3,6) \ifill f:0.7 
  \trd( 60,50,3,3,6) \ifill f:0.7 
  \rct(140,66,3,12)  \ifill f:0.7 

  \linewd 0.40 \lpatt()
  \vku(10,66,40,40,40) \vkd(10,54,40,40,40) \dl(170,66,0,-12) 
  \vku(13,63,37,40,37) \vkd(13,57,37,40,37) \dl( 13,63,0,-6)  
  \dmd(50,60,40,40)    \ifill f:0.0                           
  \dl(16,60,154,0)                                            

  \arrowheadtype t:V  \arrowheadsize l:1.2 w:1.2

  \linewd 0.30  \lpatt()
  \move(175 60) \ravec(0  -6) \ravec( 0  12) \rmove(3 -8) \htext{$4r_j$}
  \move( 18 76) \ravec(0 -12) \rmove(-2  13)              \htext{$U_j$}
  \move( 28 76) \ravec(0 -15) \rmove(-2  16)              \htext{$L_j$}
  \move( 37 77) \ravec(0 -10) \rmove(-2  11)              \htext{$S_l$}
  \move( 61 86) \ravec(0  -7) \rmove(-2   8)              \htext{$S_+$}
  \move( 61 34) \ravec(0   7) \rmove(-2 -14)              \htext{$S_-$}
  \move(142 77) \ravec(0 -10) \rmove(-2  11)              \htext{$S_r$}

  \arrowheadtype t:F  \arrowheadsize l:2.7 w:2.7

  \linewd 0.15 \lpatt(0.5 1.5)
  \dl( 10,52,0,-35) \dl( 13,56,0,-44) \dl( 16,59,0,-42) \dl(50,52,0,-40)
  \dl(130,52,0,-40) \dl(167,56,0,-39) \dl(170,52,0,-40)

  \linewd 0.30 \lpatt()
  \move(-12 8) \ravec(196   0) \rmove(-5  3) 
  \move(-12 8) \ravec(  0 100) \rmove( 2 -4) 

  \dl( 10,10, 0,-4) \rmove( -1  6) \htext{$1$}
  \dl( 13,10, 0,-4) \rmove( -6 -6) \htext{$1+r_j$}
  \dl( 16,10, 0,-4) \rmove(  0  4) \htext{$1+2r_j$}
  \dl( 50,10, 0,-4) \rmove( -6 -6) \htext{$1+R_j$}
  \dl( 90,10, 0,-4) \rmove( -6 -6) \htext{$1+2R_j$}
  \dl(130,10, 0,-4) \rmove( -6 -6) \htext{$1+3R_j$}
  \dl(167,10, 0,-4) \rmove(-27  6) \htext{$1+4R_j-r_j$}
  \dl(170,10, 0,-4) \rmove( -1 -6) \htext{$1+4R_j$}

  \dl(-14,60, 4,0 ) \rmove(  2 -2) \htext{$a_j$}
  \dl(-14,66, 4,0 ) \rmove( -1  1) \htext{$a_j+2r_j$}
  \dl(-14,54, 4,0 ) \rmove( -1 -5) \htext{$a_j-2r_j$}

 } \botcaption{\figA} The decoration $D_j$
\endcaption\medpagebreak

\example{\exaA}%
Our domain $G\subset\Re^2$ is $(0,1)^2\cup\(\bigcup_{j=1}^\infty D_j\)$,
where the sets $D_j$ are ``decorations'' attached to the right-hand side of
the unit square, centered at $(1,a_j)$. To define $D_j$, we begin with a
pair of rectangles of length $R_j$ and width $4r_j$ glued together via a
pair of bent strips of vertical width $2r_j$ that border an omitted square
of sidelength $\sqrt 2 R_j$ with sides at angle 45 degrees to the
$x_1$-axis, as in the diagram above. We then remove two closed sets. The
first removed set is the horizontal midline segment $L_j$ that begins at
the right-hand side of our decoration and ends at a distance $2r_j$ from
the left-hand side; this effectively makes the set into a union of an upper
and a lower corridor, both with three 45 degree bends. Finally, we remove a
bent U-shaped set $U_j$ which follows the (horizontal and diagonal)
midlines of the upper and lower corridors and whose points have
$x_1$-coordinates between $1+r_j$ and $1+4R_j-r_j$. Thus in our final
decoration $D_j$, there are four long bent corridors, each of vertical
width $r_j$ and with $x_1$-coordinates between $1+r_j$ and $1+4R_j$; we
call these the {\it first, second, third, and fourth corridors} in order of
increasing $y$-values. The exact values of $a_j$, $R_j$, and $r_j$ are
irrelevant as long as the decorations are pairwise disjoint and $4\le
R_j/r_j\to\infty$ as $j\to\infty$; we could for instance pick $a_j=2^{-j}$,
$R_j=4^{-j-1}$, and $r_j=8^{-j-1}$.
\endexample

The proof that $G$ is a $(0,10)$-\Wsli/ domain is a rather lengthy case
analysis similar to those in \cite{BS3, Section 3}, \cite{B1, Section 4.7},
and \cite{B2, Theorem~3.6}, so we merely mention the distinctive features of
the proof. The most interesting pairs of points $y=(y_1,y_2)$ and
$z=(z_1,z_2)$ are those that are fully contained in a single decoration
$D_j$, and which do not lie close to the boundary of the domain in the sense
that $\d_G(y),\d_G(z)\ge r_j/4$. For such points, the slices we use are one
of the following four types of slices that we collectively refer to as {\it
corridor slices}. Letting $N_j=\lfloor R_j/r_j \rfloor$, we split the part
of $D_j$ between the coordinate values $x_1=1+2r_j$ and $x_1=1+R_j$ into
$N_j$ {\it left slices} like $S_l$ all of equal width in the first
coordinate (see Figure 3.1). We similarly split the part of $D_j$ between
$x_1=1+R_j$ and $x_1=1+3R_j$ into $2N_j$ {\it upper slices} like $S_+$, and
$2N_j$ {\it lower slices} like $S_-$. Finally, we similarly split the part
of $D_j$ between $x_1=1+3R_j$ and $x_1=1+4R_j-r_j$ into $N_j$ {\it right
slices} like $S_r$.

A $(0,10)$-\Wsli/ inequality trivially holds when $k(y,z)\le 20$, so we may
assume that $k(y,z)>20$. Suppose $y,z$ both lie in the fourth corridor, and
by symmetry we assume that $y_1\le z_1$. Then we take as our admissible set
all left, upper, and right slices that lie in the set
$[y_1+r_j,z_1-r_j]\times\Re$, accompanied by their diameters. Note that
since $k(y,z)\ge 20$ and $\d_G(y),\d_G(z)\ge r_j/4$, it follows that
$z_1\ge y_1+9r_j$, and it is readily verified that the chosen slices form
a $(0,10)$ \Wsli/ dataset. As a hint note that the horizontal line segment
that runs through a left or right slice along the middle of a corridor has
quasihyperbolic length $1$. The same argument works for the other
corridors, except that we use lower slices in place of upper slices for the
first and second corridors.

If $y,z$ lie in the third and fourth corridors, respectively, we similarly
get a $(0,10)$ \Wsli/ dataset by taking all right slices that lie in the
set $[y_1\mino z_1, \infty)\times\Re$ and that do not contain points within
a distance $r_j$ of $y$ or $z$. If $y,z$ lie in the second and fourth
corridors, respectively, then we know that $k(y,z)\approx N_j$, and we get
a $(0,10)$-\Wsli/ dataset by taking all left and all right slices that do
not contain points within a distance of $r_j$ of $y$ or $z$. All
other possibilities are like one or the other of these last two cases.

Note that some or all upper and lower slices can be added to the \Wsli/
dataset for certain choices of pairs $y,z$, but not if $y,z$ are positioned
badly. For instance if $y,z$ lie in the third and fourth corridors,
respectively, and $y_1\maxo z_1\le 1+R_j$, then for every upper or lower
slice $S$, there is a path from $y$ to $z$ that avoids $S$. This problem
with ``slicing up'' the middle part of $D_j$ is precisely what makes every
slice condition fail, an argument that we now make more precise.

Suppose $\{S_i\}_{i=0}^m$ is a set of $C$-slices for the pair of points
$z:=(1+R_j,a_j+r_j/2)$, $y:=(1+R_j,a_j+3r_j/2)$, with $\g\in\G_G(y_j,z_j)$
being the associated path. Then $\g$ has to contain a point $u$ with first
coordinate $1+2R_j$ in either the first or fourth corridor. Since
$k(u,\{y,z\})$ tends to infinity as $j$ tends to infinity,
$u\in\bigcup_{i=1}^m \cl[S_i]$ if $j$ is sufficiently large. Suppose
therefore that $u\in\cl[S_i]$. Since there is a path from $y$ to $z$ that
stays a distance greater than $R_j$ from $u$, it follows from (WS-1) that
$d_i>R_j$. The slice property now ensures that $\d_G(u)\ge d_i/C>R_j/C$,
contradicting the fact that $\d_G(u)<r_j$ when $j$ is sufficiently large.
Thus $G$ is not a slice domain. \medskip

\def\open#1#2\endopen{{\noindent\bf #1. }#2}
\open{Open Problem A} Find a domain $\Om\subsetneq\Rn$ which is an inner
$0$-\wsli/ domain, but not a slice domain. More generally, one could ask
for any example of a length space which is an inner $0$-\wsli/ space but
not a slice space.
\endopen

The above problem is posed because the authors feel that slice-type
conditions, and the relationships between them, are more subtle when the
underlying metric is a length metric. Note that the previous example does
not work since the corridor slices almost all have inner Euclidean diameter
much larger than their Euclidean diameter, and so inequality (WS-1) of the
inner $(0,C)$-\wsli/ fails when $j$ is sufficiently large.

\head\AvsBsec. $\b$-\wsli/ but not $\a$-\wsli/ \endhead

\def\figB{Figure~\AvsBsec.1}
\def\figC{Figure~\AvsBsec.2}
\def\exaB{Theorem~\AvsBsec.3}
\def\exaC{Example~\AvsBsec.4}
\def\exaD{Example~\AvsBsec.5}
\def\exaE{Example~\AvsBsec.6}

Suppose $0\le\a<\b<1$. Theorem 4.1 of \cite{BS2} tells us that for domains
of product type, the inner $\b$-\Wsli/ property is equivalent to the
so-called inner $\b$-\mcig/ property. By taking the product of an interval
with Lappalainen's rather complicated examples of domains that are
$\b$-\mcig/ but not $\a$-\mcig/ \cite{L, 6.7}, we therefore get domains that
are (inner) $\b$-\Wsli/, but not (inner) $\a$-\Wsli/, whenever
$0\le\a<\b<1$. In Open Problem~B of \cite{BS2, Section 6}, the authors ask
if an $\a$-\Wsli/ domain must necessarily be a $\b$-\Wsli/ domain if
$0\le\a<\b<1$.

In this section, we answer this open problem by means of a counterexample
similar to \exaA. Another variation of this construction will give a domain
that is $\b$-\Wsli/, but not $\a$-\Wsli/, and is much simpler than the
product type domains mentioned above.

Our first two counterexamples have the form
$G:=(0,1)^2\cup\(\bigcup_{j=2}^\infty D_j\)$, where each attached decoration
$D_j$ is similar to the ones in \exaA, the only essential difference being
that the horizontal rectangular parts of $D_j$ are of width $4r_j'$ and
length $R_j'$. These altered parts are either longer and fatter, or shorter
and thinner, than before, while the diagonal parts have the same dimensions
as before. The wider corridors are pinched using linear interpolation near
where they meet narrower corridors. The following pair of diagrams of the
leftmost part of $D_j$ should suffice to make more precise what we mean.

\bigskip\medskip
\centertexdraw{
  \setunitscale 0.70

  \linewd 0.40 \lpatt()
  \move( 30 80) \lvec(150 80) \lvec(160 66) \lvec(166 72)
  \move( 30 40) \lvec(150 40) \lvec(160 54) \lvec(166 48)
  \move( 50 60) \lvec(160 60)

  \move(166 51) \lvec(160 57) \lvec(150 50) \lvec( 40 50)
  \lvec( 40 70) \lvec(150 70) \lvec(160 63) \lvec(166 69)

  \move(160 60) \lvec(166 66) \lvec(166 54) \ifill f:0.0

  \arrowheadtype t:V  \arrowheadsize l:1.2 w:1.2
  \linewd 0.30  \lpatt()
  \move(170 69) \ravec( 0 -3)   \ravec( 0 6)   \rmove(3 -6) \htext{$2r_j$}
  \move(170 51) \ravec( 0 -3)   \ravec( 0 6)   \rmove(3 -6) \htext{$2r_j$}

  \move(31 68)  \ravec(7  -8)   \move(26 69)                \htext{$U_j$}
  \move(50 72)  \ravec(8 -10)   \move(46 73)                \htext{$L_j$}

  \linewd 0.15 \lpatt(0.5 1.5)
  \dl(30,38,0,-20) \dl(40,48,0,-36) \dl(50,58,0,-40) \dl(150,38,0,-20)
  \dl(160,51,0,-39)

  \arrowheadtype t:F  \arrowheadsize l:2.7 w:2.7

  \linewd 0.30  \lpatt()
  \move(  6 8)  \ravec(172   0) \rmove(-5  3) 
  \move(  6 8)  \ravec(  0 100) \rmove( 2 -4) 

  \dl( 30,10, 0,-4) \rmove( -1  6) \htext{$1$}
  \dl( 40,10, 0,-4) \rmove( -6 -6) \htext{$1+r_j'$}
  \dl( 50,10, 0,-4) \rmove( -1  5) \htext{$1+2r_j'$}
  \dl(150,10, 0,-4) \rmove(-18  6) \htext{$1+R_j'-r_j'$}
  \dl(160,10, 0,-4) \rmove( -6 -6) \htext{$1+R_j'$}

  \dl(  4,60, 4,0 ) \rmove(  2 -2) \htext{$a_j$}
  \dl(  4,80, 4,0 ) \rmove( -1  1) \htext{$a_j+2r_j'$}
  \dl(  4,40, 4,0 ) \rmove( -1 -6) \htext{$a_j-2r_j'$}

 } \botcaption{\figB}
     The left part of $D_j$ when $R_j'>R_j$ and $r_j'>r_j$.
\endcaption\medpagebreak

\bigskip\medskip
\centertexdraw{
  \setunitscale 0.70

  \linewd 0.40 \lpatt()
  \move( 30 66) \lvec(126 66) \lvec(136 90) \lvec(146 100)
  \move( 30 54) \lvec(126 54) \lvec(136 30) \lvec(146  20)
  \move( 36 60) \lvec(126 60)

  \move(146 30) \lvec(136 40) \lvec(126 57) \lvec( 33 57)
  \lvec( 33 63) \lvec(126 63) \lvec(136 80) \lvec(146 90)

  \move(126 60) \lvec(136 70) \lvec(146 80) \lvec(146 40) \lvec(136 50)
  \ifill f:0.0

  \arrowheadtype t:V  \arrowheadsize l:1.2 w:1.2
  \linewd 0.30  \lpatt()
  \move(150 90) \ravec( 0 -10)   \ravec( 0 20)  \rmove(3 -12) \htext{$2r_j$}
  \move(150 30) \ravec( 0 -10)   \ravec( 0 20)  \rmove(3 -12) \htext{$2r_j$}

  \move(40 71)  \ravec( 6  -7)   \move(36 72)                 \htext{$U_j$}
  \move(50 71)  \ravec( 8 -10)   \move(46 72)                 \htext{$L_j$}

  \linewd 0.15 \lpatt(0.5 1.5)
  \dl(30,52,0,-34) \dl(33,55,0,-43) \dl(36,58,0,-40) \dl(126,52,0,-34)
  \dl(136,28,0,-16)

  \arrowheadtype t:F  \arrowheadsize l:2.7 w:2.7

  \linewd 0.30 \lpatt()
  \move( 4 8) \ravec(156   0) \rmove(-2  3) 
  \move( 4 8) \ravec(  0 100) \rmove( 2 -4) 

  \dl( 30,10, 0,-4) \rmove(-1  6) \htext{$1$}
  \dl( 33,10, 0,-4) \rmove(-6 -6) \htext{$1+r_j'$}
  \dl( 36,10, 0,-4) \rmove(-1  5) \htext{$1+2r_j'$}
  \dl(126,10, 0,-4) \rmove(-7  6) \htext{$1+R_j'$}
  \dl(136,10, 0,-4) \rmove(-6 -6) \htext{$1+R_j'+r_j'$}

  \dl(  4,60, 4,0 ) \rmove( 2 -2) \htext{$a_j$}
  \dl(  4,66, 4,0 ) \rmove(-1  1) \htext{$a_j+2r_j'$}
  \dl(  4,54, 4,0 ) \rmove(-1 -6) \htext{$a_j-2r_j'$}

 } \botcaption{\figC}
     The left part of $D_j$ when $R_j'<R_j$ and $r_j'<r_j$.
\endcaption\medpagebreak

Let us take $R_j:=2^{-jp}$, $R_j':=2^{-jp'}$, $r_j:=2^{-jq}$,
$r_j':=2^{-jq'}$, where the quadruple $(p,q,p',q')$ is {\it allowable} if
$0<p\le q-2$, $0<p'\le q'-2$, $p\ge 2$, and $q'\ge 2$; the last two bounds
are assumed merely to ensure that we can attach all these decorations to
one side of the unit square without overlap. The exact locations of the
decorations, i.e.~the values of $a_j$, are irrelevant as long as they do
not overlap.

\proclaim{\exaB} Given $0<\a_0<1$, any allowable choice of $p,p',q,q'$ with
$p'=p+1-\a_0$ and $q'=q+1$ gives a domain $G$ which is an $\a$-\Wsli/
domain for $\a\le\a_0$, but not an $\a$-\wsli/ domain for $\a>\a_0$.
\endproclaim

\demo{Sketch of proof}%
Writing $N_j'=\lfloor R_j'/r_j'\rfloor$ and $N_j=\lfloor R_j/r_j\rfloor$,
we define corridor slices as in \exaA, so that there are $N_j'$ left slices
between $x_1=1+2r_j'$ and $x_1=1+R_j'-r_j'$, $2N_j$ upper and $2N_j$ lower
slices between $x_1=1+R_j'$ and $x_1=1+R_j'+2R_j$, and $N_j'$ right slices
between $x_1=1+R_j'+2R_j+r_j'$ and $x_1=1+2R_j'+2R_j-r_j'$. In this proof,
$A\ll B$ means that $A/B\to 0$ as $j\to\infty$.

Note that the $d_{\a,G}$-length of the (horizontal or diagonal) line
segment given by the intersection of a single corridor slice with the
midline of that corridor is comparable with $r_j^\a$ for an upper or lower
slice $S$, and $(r_j')^\a$ for a left or right slice $S$. For some pairs of
points $y,z$, the $(\a,C)$-\Wsli/ defining inequality holds using a similar
argument to that in \exaA\ once we pick $C=C(\a)$ to be large enough.
However this method fails in other cases. The basic obstacle is revealed by
taking $y=(1+R_j'/2,a_j+3r_j'/2)$ and $z=(1+R_j'/2,a_j+r_j'/2)$. Then all
paths from $y$ to $z$ have to go through complete horizontal and diagonal
parts of at least two corridors, and so it follows that $d_\a(y,z)\approx
L_{j,\a}+ L_{j,\a}'$, where $L_{j,\a}:=R_j/r_j^{1-\a}$ and
$L_{j,\a}':=R_j'/(r_j')^{1-\a}$. We cannot use upper or lower slices in any
admissible set for $y,z$ since there always exist connecting paths that
avoid any given set of this type, but the set of all right slices $S$
(paired with their diameters $d_S$) always gives a $10$-admissible set.
Denoting by $\F$ the set of such right slices, we see that $\sum_{S\in\F}
d_S^\a\approx L_{j,\a}'$. Since
$$
{1\over j}\cdot\log_2\({L_{j,\a}\over L_{j,\a}'}\) = p'-p-(1-\a)(q'-q) =
\a-\a_0,
$$
we see that $\W_0:=\{(S,\diam(S)\mid S\in\F\}$ is an $(\a,C)$-\Wsli/
dataset for appropriate $C=C(\a)$ as long as $\a\le\a_0$, as required.
However, $\W_0$ fails to be an $(\a,C)$-\Wsli/ dataset when $\a>\a_0$ since
then $L_{j,\a}'\ll L_{j,\a}$.

Given $\a\in(\a_0,1)$, it remains to show that there are no $(\a,C)$-\wsli/
datasets for the pair $y,z$, assuming that $j$ is sufficiently large.
Suppose for the sake of contradiction that $\W:=\{(S,d_S)\mid S\in\F\}$ is
some $(\a,C)$-\wsli/ dataset and write $\Sum_\a:=\sum_{S\in\F} d_S^\a$.
Since $\d^\a(y)+\d^\a(z)\approx r_j^\a\ll L_{j,\a}$, it follows that
$\Sum_\a\approx d_\a(y,z)\approx L_{j,\a}$. Using (WS-1) and the geometry
of the domain, we see that any slice that includes points outside $D_j$ has
diameter larger than $R_j'/2$. Furthermore, if $m_i$ is the number of such
slices $S$ for which $d_S\in(2^{i-1}R_j',2^i R_j']$, then (WS-1) and the
fact that there are paths from $y$ to $z$ of length comparable to $R_j$
together imply that $m_i\lc 2^{-i}R_j/R_j'$. By summing the resulting
series over the index $i$, we see that the contribution of all such slices
to $\Sum_\a$ is at most comparable with $A_j:=R_j/(R_j')^{1-\a}$, and
$A_j\ll L_{j,\a}$ because $p'=p+1-\a_0<q$. We can therefore delete these
slices from our dataset and our redefined set $\W$ is still an
$(\a,C)$-\wsli/ dataset (if we suitably redefine $C$).



Consider next from the remaining slices those that do not enter into any
diagonal corridor by a distance more than $R_j'$ from the base (meaning the
left and right ends of the diagonal corridors of $D_j$). We let $\lambda$
temporarily denote the path in $\Gamma_G(y,z)$ that runs along the U-shaped
mid-corridor path on the right. Since the intersection of $\lambda$ with the
slices under present consideration can have length at most comparable to
$R_j'$, it follows that the number $m_i$ of such slices $S$ for which
$d_S\in(2^i r_j',2^{i+1} r_j']$ is at most comparable to $2^{-i}R_j'/r_j'$.
Since any such slice has diameter at least comparable to $r_j'$, it follows
that the contribution of such slices is at most comparable to $L_{j,\a}'\ll
L_{j,\a}$.

It remains to consider the slices which lie in $D_j$, and enter into at
least one diagonal corridor by a distance exceeding $R_j'$ from the base.
Let $m_i$ be the number of such slices $S$ for which $d_S\in(2^i
R_j',2^{i+1}R_j']$. Now such slices must include points that are a distance
at most comparable with $2^iR_j'$ from the base, since if all points in such
a slice are much further than this from the base then the slice cannot
contain points in both the upper and lower pair of corridors and so cannot
separate the pair $y,z$, contradicting (WS-1). We deduce that such slices
are fully contained within a distance comparable to $2^iR_j'$ of the base,
and so the intersection of $\lambda$ with such slices can have length at
most comparable to $2^i R_j'$. It follows that $m_i\lc 1$. In order to
accommodate all such slices, the index $i$ need only run up to the value
$\log_2(R_j/R_j')$. Consequently, we may estimate the contributions of these
remaining slices with the upper bound: $\sum_{i=0}^{\log_2(R_j/R_j')}1\cdot
(2^iR_j')^\alpha \approx (R_j')^\alpha\sum_{i=0}^{\log_2(R_j/R_j')}2^{\alpha
i}\approx (R_j')^\alpha(R_j/R_j')^\alpha \approx R_j^\alpha$.   But this is
much smaller than $L_{j,\a}$ when $j$ is large and so we get a
contradiction. \QED\enddemo

The above construction is quite flexible: it can be varied to give examples
with various other types of behavior. We content ourselves below with three
variants, but first let us define $\a(D)$, the {\it $\a$-set of a domain
$D\subsetneq\Rn$}, to be the set of all $\a\in[0,1)$ for which a given
domain $D$ is an $\a$-\wsli/ domain. \exaB\ shows that there are domains
$G$ with $\a(G)=[0,\a_0]$ for each $0<\a_0<1$. By varying some of the
details in the definition of $G$, we now get some other $\a$-sets. We omit
the details of the proofs which are all similar to that of \exaB.

The first of our three examples allows us to get the same $\a$-sets as the
product-type examples mentioned at the beginning of this section, but is
much simpler. Our second example shows that the endpoint of our $\a$-set can
be omitted, and the third shows that $\a$-sets need not be intervals.

\proclaim{\exaC} If $0<\a_0<1$, then any allowable choice of $p,p',q,q'$
with $p=p'+1-\a_0$ and $q=q'+1$ gives a domain $G$ with $\a(G)=[\a_0,1)$.
\endproclaim

\example{\exaD} If we redefine $r_j:=j2^{-jq}$ in \exaB, but leave
everything else unchanged, then $\a(G)=[0,\a_0)$. The key fact is that when
$\a=\a_0$, we now have $L_{j,\a}/L_{j,\a}'\to\infty$ as $j\to\infty$.
\endexample

\example{\exaE} Consider a domain with decorations $D_j$ similar to those of
\exaB, but with two rectangular parts on both sides of the diagonal part.
The diagonal part and the innermost pair of rectangular parts of $D_j$ are
identical in shape to the full decoration $D_j$ of \exaB\ with the exception
that we must alter $U_j$ and $L_j$ so that they also pass through the outer
rectangular parts, which have length $R_j'':=2^{-jp''}$ and width
$4r_j'':=4\cdot2^{-jq''}$. These outer parts are chosen to be longer and
fatter than the inner rectangular parts and are connected by linear
interpolation to the inner parts as before. By choosing $p'=p+1-\a_0$ and
$q'=q+1$, and $p''=p-1+\a_1$ and $q''=q-1$, $0<\a_0<\a_1<1$, it follows that
$\a(G)=[0,\a_0] \cup [\a_1,1)$.
\endexample


By taking $(p,q)=(3,6)$ and $(p',q')=(p,q)\pm(1-\a_0,1)$ for some
$\a_0\in[0,1)$, it is clear that $(p,q,p',q')$ is always allowable. This
allows us to consider domains consisting of a sequence of decorations $D_j$
joined to the unit square that generalize the above constructions. Each
$D_j$ has diagonal corridors specified by the dimensional parameters
$R_j:=2^{-3j}$ and $r_j:=2^{-6j}$, and $D_j$ also has one or more horizontal
corridors on each side of these diagonal corridors, symmetrically
distributed around the center of the diagonal corridor: if the $i$th
horizontal corridor on the left counting outwards from the diagonal corridor
has dimensional parameters $R_{j,i}:=2^{-jp_{j,i}}$ and
$r_{j,i}:=2^{-jq_{j,i}}$, then the $i$th horizontal corridor on the right is
defined by these same parameters. Here $(p_{j,i},q_{j,i})-(3,6)$ is always
$\pm(1-\a_{j,i},1)$ for some $0<\a_{j,i}<1$. The corridors are joined by
linear interpolation as before. We call these {\it corridor decorations} and
we call the domain $\Om$ obtained by joining a sequence of such corridor
decorations a {\it decorated square (with corridor decorations
$(D_j)_{j=1}^\infty$)}.

It is not hard to show that the sets $\a(\Om)$ for the set of decorated
squares $\Om$, are closed under countable intersections and finite unions.
Let us justify this first for intersections. Suppose $(\Om_k)_{k=1}^\infty$
is a sequence of decorated squares with corridor decorations
$(D_{k,j})_{j=1}^\infty$. It is routine to show that we can define a
decorated square $\Om$ with corridor decorations $(D_i)_{i=1}^\infty$, where
$D_i:=D_{k_i,j_i}$ for some appropriate choice of $k_i,j_i$, such that
$\a(\Om)=\bigcap_{k=1}^\infty \a(\Om_k)$.

As for finite unions, if we have a finite set of decorated squares $\Om_k$,
$k=1,\dots,k_0$, with corridor decorations $(D_{k,j})_{j=1}^\infty$, then we
take our cue from \exaE: for fixed $j$, we join together the horizontal
corridors of each $D_{k,j}$, $k=1,\dots,k_0$, as we did in \exaE\ to get a
new decoration $D_j$. The decorated square $\Om$ with corridor decorations
$(D_j)_{j=1}^\infty$ then has the property that $\a(\Om)=\bigcup_{k=1}^{k_0}
\a(\Om_k)$.

The above constructions suggest that every Borel subset of $[0,1)$ may well
be of the form $\a(G)$ for some bounded domain $G\subset\Rn$. However we do
not know if this is so.

As pointed out at the start of this section, there are domains in $\Rn$
which are inner $\b$-\Wsli/ but not inner $\a$-\Wsli/ whenever
$0\le\a<\b<1$. However none of our decorated examples above are inner
$\a$-\wsli/ domains, so they cannot answer the following problem.

\medskip
\open{Open Problem B} Given $0\le\a<\b<1$, is there a domain in $\Rn$ which
is inner $\a$-\Wsli/ (or even just $\a$-\wsli/) but not $\b$-\Wsli/? More
generally, one could ask for any example of a length space which is inner
$\a$-\Wsli/ (or even just $\a$-\wsli/) but not $\b$-\Wsli/.
\endopen


\enddocument

\Refs
\widestnumber\key{BK2}

\ref \key BB \by Z. Balogh and S.M. Buckley \paper Geometric
characterizations of Gromov hyperbolicity \jour Invent. Math. \vol 153 \yr
2003 \pages 261--301 \endref

\ref \key BHK \by M. Bonk, J. Heinonen, and P. Koskela \paper Uniformizing
Gromov hyperbolic spaces \jour Ast\'erisque \vol 270 \endref

\ref \key B1 \by S.M. Buckley \paper Slice conditions and their applications
\jour Rep. Univ. Jyv\"askyl\"a Dept. Math. Stat. \vol 92 \yr 2003 \pages
63--76 \endref

\ref \key B2 \by S.M. Buckley \paper Quasiconformal images of H\"older
domains \jour Ann. Acad. Sci. Fenn. \vol 29 \yr 2004 \pages 21--42 \endref

\ref \key BK1 \by S.M. Buckley and P. Koskela \paper Sobolev-Poincar\'e
implies John \jour Math. Research Letters \vol 2 \pages 577--593 \yr
1995 \endref

\ref \key BK2 \by S.M. Buckley and P. Koskela \paper Criteria for
Imbeddings of Sobolev-Poincar\'e type \jour Internat. Math. Res. Notices
\yr 1996 \pages 881--901 \endref

\ref \key BO \by S.M. Buckley and J. O'Shea \paper Weighted
Trudinger-type inequalities \jour Indiana Univ. Math. J. \vol 48 \yr 1999
\pages 85--114 \endref

\ref \key BS1 \by S.M. Buckley and A. Stanoyevitch \paper Weak slice
conditions and H\"older imbeddings \jour J. London
Math. Soc. \vol 66(3) \pages 690--706 \yr 2001 \endref

\ref \key BS2 \by S.M. Buckley and A. Stanoyevitch \paper Weak slice
conditions, product domains and quasiconformal mappings \jour
Rev. Mat. Iberoamericana \vol 17 \yr 2001 \pages 1--37 \endref

\ref \key BS3 \by S.M. Buckley and A. Stanoyevitch \paper Distinguishing
properties of weak slice conditions \jour Conform. Geom. Dyn. \vol 7 \yr
2003 \pages 49--75 \endref

\ref \key GM \by F.W. Gehring and O. Martio \paper Lipschitz classes and
quasiconformal mappings \jour Ann. Acad. Sci. Fenn. Ser. A I Math. \vol
10 \yr 1985 \pages 203--219 \endref

\ref \key GO \by F.W. Gehring and B. Osgood \paper Uniform domains and
the quasihyperbolic metric \jour J. Analyse Math. \vol 36 \yr 1979
\pages 50--74 \endref

\ref \key HK \by J. Heinonen and P. Koskela \paper Quasiconformal maps
in metric spaces with controlled geometry \jour Acta Math. \vol 181
\yr 1998 \pages 1--61 \endref

\ref \key L \by V. Lappalainen \paper $Lip_h$-extension domains
\jour Ann. Acad. Sci. Fenn. Ser. A I Math. Diss. \vol 56 \pages 1--52
\yr 1985 \endref

\ref \key Mz \by V.L. Maz'ya \book Sobolev Spaces \publ Springer-Verlag \yr
1985 \publaddr Berlin \endref

\ref \key V1 \by J. V\"ais\"al\"a \paper On the null-sets for extremal
distances \jour Ann. Acad. Sci. Fenn. Ser. A I \vol 322 \yr 1962
\pages 12pp \endref

\ref \key V2 \by J. V\"ais\"al\"a \book Lectures on $n$-dimensional
quasiconformal mappings \bookinfo Lecture Notes in Mathematics 229
\publ Springer-Verlag \publaddr Berlin \yr 1970
\endref

\ref \key V3 \by J. V\"ais\"al\"a \paper Uniform domains \jour Tohoku
Math. J. \vol 40 \yr 1988 \pages 101--118 \endref

\ref \key V4 \by J. V\"ais\"al\"a \paper Quasiconformal mappings of
cylindrical domains \jour Acta Math. \vol 162 \yr 1989 \pages 201--225
\endref

\ref \key V5 \by J. V\"ais\"al\"a \paper Relatively and inner uniform
domains \jour Conf. Geom. Dyn. \vol 2 \yr 1998 \pages 56--88
\endref

\endRefs
\bye